\documentclass[a4paper]{article}
\newcommand{\la}{\lambda}
\newcommand{\x}{\xi}
\newcommand{\h}{\eta}
\newcommand{\De}{\Delta}
\newcommand{\Ga}{\Gamma}
\newcommand{\ga}{\gamma}
\newcommand{\al}{\alpha}
\newcommand{\bV}{{\bf V}}
\newcommand{\bW}{{\bf W}}
\newcommand{\bu}{{\bf u}}
\newcommand{\bv}{{\bf v}}
\newcommand{\bx}{{\bf x}}

\newcommand{\tensor}{\otimes \!_Z}
\newcommand{\Pz}{{\cal P}_Z}
\newcommand{\Pk}{{\cal P}_K}
\newcommand{\Pl}{{\cal P}_L}
\newcommand{\Pkh}{\widehat{{\cal P}}_K}
\newcommand{\Slk}{{\cal S}_{L/K}}
\newcommand{\AkL}{{\cal A}_K(L)}
\newcommand{\cC}{{\cal C}}
\newcommand{\cZ}{{\cal Z}}
\newcommand{\cF}{{\cal F}}
\newcommand{\cFt}{\widetilde{\cal F}}
\newcommand{\cA}{{\cal A}}
\newcommand{\cAt}{\widetilde{\cal A}}
\newcommand{\cX}{{\cal X}}
\newcommand{\cQ}{{\cal Q}}
\newcommand{\cH}{{\cal H}}
\newcommand{\AGL}{\mbox{{\rm AGL}}(1,L)}
\newcommand{\quer}{(\overline{\mbox{ \rule{0ex}{1.4ex} }})}
\newcommand{\beweisende}{~\rule{1,2ex}{1,2ex}}
\newtheorem{theo}{Theorem}
\newtheorem{lemma}{Lemma}
\title{Affine Circle Geometry over\\Quaternion Skew Fields}
\author{By Hans Havlicek}
\date{}

\begin{document}

\maketitle

\begin{abstract}
We investigate the affine circle geometry arising from a quaternion
skew field and one of its maximal commutative subfields.
\end{abstract}

\thispagestyle{empty}
%%%%%%%%%%%%%%%%%%%%%%%%%%%%%%%%%%%%%%%%%%%%%%%%%%%%%%%%%%%%%%%%%%%%%%
%additionally for page 1
%%%%%%%%%%%%%%%%%%%%%%%%%%%%%%%%%%%%%%%%%%%%%%%%%%%%%%%%%%%%%%%%%%%%%%
\section{Introduction}\label{sect-INTRO}

\subsection{}
The present paper is concerned with the chain geometry $\Sigma(K,L)$
(cf. \cite{Benz1}) on a field extension $L/K$, where $K$ is a maximal
commutative subfield of a quaternion skew field $L$. Thus $L$ is not a
$K$--algebra. This has many geometric consequences. Best known is
probably that three distinct points do not determine a unique chain.
As in ordinary M\"obius--geometry, it is possible to obtain an affine
plane by deleting one point, but a more sophisticated technique is
necessary in order to define the lines of this plane. We take a closer
look on this construction from two different points of view, starting
either from a spread of lines associated to $\Sigma(K,L)$ or the point
model of this spread on the Klein quadric. The chains of $\Sigma(K,L)$
yield the lines, degenerate circles and non--degenerate circles of
such an affine plane. We establish some properties of these circles
and show that degenerate circles are affine Baer subplanes. If $K$ is
Galois over the centre of $L$ then each non--degenerate circle can be
written as intersection of two affine Hermitian varieties.

We encourage the reader to compare our results with the survey article
\cite{Herz1} on chain geometry over an algebra and \cite{Herz2}. There
is an extensive literature on the real quaternions. A lot of
references can be found, e.g., in \cite{Benz1}, \cite{Berg1},
\cite{HH4}, \cite{TaSca1}, \cite{Wilk1}.

\subsection{}
Throughout this paper $L$ will denote a quaternion skew field with
centre $Z$ and $K$ will be a maximal commutative subfield of $L$. The
following exposition follows \cite{Cohn1}, \cite{HH3},
\cite[pp.168--171]{Pick1}.

Choose any element $a\in K\setminus Z$ with minimal equation, say
   \begin{displaymath}
   a^2 + a\la_1 + \mu_1 = 0 \quad (\la_1,\mu_1 \in Z).
   \end{displaymath}

If $K/Z$ is Galois then
   \begin{displaymath}
   \quer: K \to K, \quad u=\x+a\h \mapsto \overline{u} :=
   \x-(\la_1+a)\h \quad (\x,\h \in Z)
   \end{displaymath}
is an automorphism of order 2 fixing $Z$ elementwise%
   \footnote
   {By an appropriate choice of $a$ it would be possible to have
   $\la_1=0$ (Char$K\neq 2$) or $\la_1=1$ (Char$K=2$).}%
. There exists an element
$i \in L \setminus K$ such that
   \begin{displaymath}
   i^{-1} u i = \overline{u}\quad \mbox{for all } u \in K,
   \end{displaymath}
whence
   \begin{equation}\label{galois-comm}
   u i = i \overline{u} \quad \mbox{for all } u \in K.
   \end{equation}

If $K/Z$ is not Galois then, obviously, $\mbox{Char}K=2$ and
$\la_1=0$. The mapping
   \begin{displaymath}
   D : K \to K, \quad u=\x+a\h \mapsto u^D:=a\h \quad (\x,\h \in Z)
   \end{displaymath}
is additive and satisfies $(uu')^D = u^D u' + u u'^D$ for all $u,u'
\in K$, i.e., $D$ is a derivation of $K$. There exists an $i \in
L\setminus K$ such that
   \begin{displaymath}
   a^{-1}ia = i+1
   \end{displaymath}
which leads to the rule
   \begin{equation}\label{non-galois-comm}
   ui = i u + u^D \quad \mbox{for all } u \in K.
   \end{equation}

In every case the element $i$ has a minimal equation over $Z$, say
   \begin{displaymath}
   i^2 + i\la_2 + \mu_2 = 0 \quad (\la_2,\mu_2 \in Z).
   \end{displaymath}
If $K/Z$ is Galois then $i^2 \in Z$, whence $\la_2=0$. If $K/Z$ is not
Galois then $i$ and $i+1$ have the same minimal equation. This implies
$\la_2=1$. The mapping
   \begin{equation}
   A:L\to L,\quad u+iv  \mapsto
   \left\{
      \begin{array}{l@{\mbox{ : }}l}
      \overline{u}-iv     & K/Z\mbox{ Galois,}\\
      u+v+vi              & K/Z\mbox{ not Galois,}
      \end{array}
   \right\}
   \quad (u,v\in K)
   \end{equation}
is an involutory antiautomorphism of $L$ fixing $K$. The norm of $x\in
L$ is given by $N(x):=x^Ax$.

\subsection{}
The mappings $\quer$ and $D$ allow, respectively, the following
geometric interpretations:

Let $\bV$ be a right vector space over $Z$, $\dim\bV\geq 2$. We are
extending $\bV$ to $\bV \tensor K$ with $\bv \in \bV$ to be identified
with $\bv\otimes 1$. Then define a mapping $\bV\tensor K \to
\bV\tensor K$ by
   \begin{displaymath}
   \sum_{\bv \in \bV}\bv\otimes k_\bv \mapsto
   \left\{
      \begin{array}{l@{\mbox{ : }}l}
      \sum\limits_{\bv \in \bV}\bv\otimes \overline{k}_\bv
               & K/Z\mbox{ Galois,}\\
      \sum\limits_{\bv \in \bV}\bv\otimes k_\bv^D
               & K/Z\mbox{ not Galois,}
      \end{array}
   \right\}
   \quad(k_\bv \in K).
   \end{displaymath}
By abuse of notation, this mapping will also be written as $\quer$ and
$D$, respectively.

In terms of the projective spaces $\Pz(\bV)$ and $\Pk(\bV\tensor K)$
the first projective space is being embedded in the second one as a
Baer subspace. If $\bx K$ is a point of $\Pk(\bV\tensor
K)\setminus\Pz(\bV)$ then through this point there is a unique line of
$\Pk(\bV\tensor K)$ containing more than one point of $\Pz(\bV)$. That
line is given by
   \begin{displaymath}
   \bx K\vee \overline{\bx}K \quad \mbox{and} \quad \bx K\vee (\bx^D)
   K,
   \end{displaymath}
respectively%
   \footnote
   {At least in the first case this is very well known.}%
. Note that defining a mapping by setting $\bx K\mapsto (\bx^D)K$ is
ambiguous, since
   \begin{displaymath}
   (\bx u)^D = \bx^D u + \bx u^D \quad\mbox{for all } \bx\in\bV\tensor
   K\mbox{, }u\in K.
   \end{displaymath}

We give a second interpretation in terms of affine planes%
   \footnote
   {Cf. the concept of `Minimalkoordinaten' described, e.g., in
   \cite[p.35]{Wund1}}%
:

   \begin{lemma}\label{lemma-b}
   Let $\bW$ be a right vector space over $K$, $\dim\bW=2$, and let
   $\{\bu,\bv\}$ be a basis of $\bW$. Then
      \begin{equation}\label{minimal}
         \begin{array}{l@{\mbox{ : }}l}
         \{ \bu k + \bv \overline{k} \mid k\in K \} &  K/Z\mbox{
         Galois,}\\
         \{ \bu k + \bv k^D\mid k\in K \}           &  K/Z\mbox{ not
         Galois,}
         \end{array}
      \end{equation}
   is an affine Baer subplane (over $Z$) of the affine plane on $\bW$.
     \end{lemma}
{\it Proof}. If $\bu'$ and $\bv'$ are linearly independent vectors of
$\bW$ then the set of all linear combinations of $\bu'$ and $\bv'$
with coefficients in $Z$ is an affine Baer subplane over $Z$. Write
$k=\x+a\h$ with $\x,\h\in Z$.

If $K/Z$ is Galois then
   \begin{displaymath}
   \bu k + \bv \overline{k} = (\bu+\bv)\xi + (\bu a - \bv(\la_1 +
   a))\h.
   \end{displaymath}
The vectors $\bu+\bv$ and $\bu a - \bv(\la_1 + a)$ are linearly
independent, since otherwise we would have the contradiction
$\overline{a} = -\la_1 - a = a$.

If $K/Z$ is not Galois then
   \begin{displaymath}
   \bu k + \bv k^D = \bu\xi + (\bu + \bv)a\h.
   \end{displaymath}
The vectors $\bu$ and $(\bu + \bv)a$ are linearly
independent.\beweisende

\section{Projective Chain Geometry on $L/K$}\label{sect-B}

\subsection{}
Let $L/K$ be given as before. Following \cite[p.320ff.]{Benz1} we
obtain an incidence structure $\Sigma(K,L)$ as follows: The points of
$\Sigma(K,L)$ are the points of the projective line over $L$, viz.
$\Pl (L^2)$, the blocks, now called chains, are the $K$-sublines of
$\Pl (L^2)$. However, in contrast to \cite{Benz1}, we shall regard
$L^2$ as right vector space over $K$ rather than $L$. Each
$s$--dimensional subspace of $L^2$ (over $L$) is $2s$--dimensional
over $K$, whence $\Pk(L^2)=:\Pk$ is 3--dimensional. The points of
$\Sigma(K,L)$ now appear as lines of a spread of $\Pk$, say $\Slk$;
cf. \cite{HH1}, \cite{HH2}. If $t$ is a line of $\Pk$ not contained in
$\Slk$ then through each point of $t$ there goes exactly one line of
$\Slk$. The subset $\cC$ of $\Slk$ arising in this way is a chain of
$\Sigma(K,L)$. We call $t$ a transversal line of the chain $\cC$.  If
$L/K$ is not Galois then each chain has exactly one transversal line,
otherwise exactly two transversal lines that are interchanged under
the non--projective collineation
   \begin{equation}\label{def-iota}
   \iota : \Pk \to \Pk,\quad (l_0,l_1)K \mapsto (l_0 i,l_1 i)K.
   \end{equation}
Cf. \cite[Theorem 2]{HH3}, \cite{KiRei1}.

     \subsection{}
Write $\cal L$ for the set of lines of $\Pk$ and $\ga : {\cal L} \to
\Pkh$ for the Klein mapping. Here $\Pkh$ is the ambient space of the
Klein quadric $\cQ:={\cal L}^\ga$. The underlying vector space of
$\Pkh$ is $L^2\wedge L^2$ (over $K$). In \cite[Theorem 1]{HH3} it is
shown that there is a unique 5--dimensional Baer subspace $\Pi_Z$
(over $Z$) of $\Pkh$ such that
   \begin{displaymath}
   \Slk{}^{\ga}=\Pi_Z\cap{\cQ}.
   \end{displaymath}
With respect to $\Pi_Z$ the set $\Slk{}^{\ga}$ is an oval quadric,
i.e. a quadric without lines. A subset $\cC$ of $\Slk$ is a chain if,
and only if, there exists a 3--dimensional subspace $\cX$ of $\Pkh$
such that%
   \footnote
{If $L$ is the skew field of real quaternions then $K$ is a field of
complex numbers and $Z$ the field of real numbers. Here conditions
(\ref{sphere-1}) and (\ref{sphere-2}) are already sufficient to
characterize the $\ga$--images of chains.}
   \begin{eqnarray}
   \label{sphere-1}
   {} & {}& {\cX}\cap \Pi_Z \mbox{ is a 3--dimensional subspace of
   }\Pi_Z,\\
   \label{sphere-2}
   {} & {}& {\cC^{\ga}}={\cX}\cap \Slk{}^{\ga} \mbox { is an elliptic
   quadric of }{\cX}\cap\Pi_Z \mbox{ (over }Z),\\
   \label{sphere-3}
   {} & {} & {\cX}\cap{\cQ} \mbox{ contains a line of }\Pkh;
   \end{eqnarray}
cf. \cite[Theorem 1]{HH3}.

\subsection{}
The automorphism group of $\Sigma(K,L)$ is formed by all bijections of
$\Slk$ taking chains to chains in both directions. If $\kappa$ is a
collineation or a duality of $\Pk$ with $\Slk{}^{\kappa}=\Slk$ then
$\kappa$ is yielding an automorphism of $\Sigma(K,L)$. Conversely,
according to \cite{MMN1} and \cite[Theorem 4]{HH3}, each automorphism
of $\Sigma(K,L)$ can be induced by an automorphic collineation or
duality of $\Slk$, say $\kappa$. This $\kappa$ is uniquely determined
for $K/Z$ not being Galois, otherwise the product of $\iota$ (cf.
formula (\ref{def-iota})) and $\kappa$ is the only other solution.

Transferring these results to $\Pkh$ establishes that an automorphic
collineation $\mu$ of the Klein quadric is the $\ga$--transform of an
automorphic collineation or duality of $\Slk$ if, and only if, $\Pi_Z$
is invariant under $\mu$. If $K/Z$ is Galois, then the
$\ga$--transform of the collineation $\iota$ (cf. (\ref{def-iota})) is
the Baer involution of $\Pkh$ fixing $\Pi_Z$ pointwise. See
\cite[Theorem 4]{HH3}.

\subsection{}
Let $\cC_0$ and $\cC_1$ be two chains with a common element, say
$p\in\Slk$. We say that $\cC_0$ is {\bf tangent} to $\cC_1$ at $p$ if
there exist transversal lines $t_i$ of $\cC_i$ $(i=0,1)$ such that
$p,t_0,t_1$ are in one pencil of lines. This is a reflexive and
symmetric relation.

If $K/Z$ is Galois then there is also an orthogonality relation on the
set of chains: If $\cC_i$ ($i=0,1$) are chains with transversal lines
$t_i$, $t_i{}^{\iota}$, respectively, then $\cC_0$ is said to be {\bf
orthogonal} to $\cC_1$ if $t_0$ intersects both $t_1$ and
$t_1{}^{\iota}$. This relation is symmetric, since $\iota$ is an
involution. Given two orthogonal chains their transversal lines form a
skew quadrilateral.

The two definitions above are not given in an intrinsic way. However,
both relations are invariant under automorphic collineations and
dualities of $\Slk$ and hence invariant under automorphisms of
$\Sigma(K,L)$.

The proofs of the following results are left to the reader: Chains
$\cC_0$, $\cC_1$ are tangent at $p\in\cC_0\cap\cC_1$ if, and only if,
their images under the Klein mapping are quadrics with the same
tangent plane at the point $p^\ga$. A chain $\cC_0$ is orthogonal to a
chain $\cC_1$ if, and only if, the subspace of $\Pkh$ spanned by
$\cC_0{}^\ga$ contains the orthogonal subspace (with respect to the
Klein quadric) of $\cC_1{}^\ga$.

\section{Affine Circle Geometry on $L/K$}\label{sect-C}

\subsection{}
With the notations introduced in section \ref{sect-B}, select one line
of $\Slk$ and label it $\infty$. Let $\cAt$ be a (projective) plane of
$\Pk$ through $\infty$ and write $\cA:=\cAt\setminus \infty$. Then
$\cA$ can be viewed as an affine plane with $\infty$ as line at
infinity. The mapping
   \begin{equation}\label{def-rho}
   \rho : \Slk\setminus \{ \infty \} \to \cA,
   \quad s \mapsto \cA\cap s
   \end{equation}
is well--defined and bijective. A chain $\cC$ containing $\infty$
yields an affine line $(\cC\setminus \{ \infty \})^{\rho}$ if, and
only if, ${\cC}$ has a transversal line in $\cAt$. Two chains with
transversal lines in $\cAt$ yield parallel lines if, and only if, the
chains are tangent at $\infty$.

If $\cAt'$ is any plane through $\infty$ then, with $\cA' :=
\cAt'\setminus \infty$, the mapping
   \begin{displaymath}
   \beta : \cA\to {\cA}',\quad \cA\cap s\mapsto {\cA}'\cap s
   \quad (s \in\Slk\setminus \{ \infty \})
   \end{displaymath}
is a well--defined bijection%
   \footnote
   {One could also select some point $A\in \infty$ and then obtain an
   affine plane by a dual construction.}%
. This $\beta$ is an affinity if either $\cAt' = \cAt$ or $\cAt' =
\cAt^{\iota}$ \cite[Theorem 5]{HH1}; the second alternative is only
possible when $K/Z$ is Galois.

\subsection{}
The group of automorphic collineations of $\Slk$ operates 3--fold
transitively on the lines of $\Slk$ \cite[p.322]{Benz1}. Thus we may
transfer $\infty$ to the line given by $(0,1)L$. Moreover, for all
$c\in L$, $c\neq 0$
   \begin{displaymath}\label{c-mal}
   (l_0,l_1)K \mapsto (c l_0,c l_1)K \quad ((0,0)\neq (l_0,l_1) \in
   L^2)
   \end{displaymath}
is an automorphic collineation of $\Slk$ fixing $\infty$. Hence,
without loss of generality, we may assume in the sequel that
   \begin{displaymath}
   \infty=\Pk((0,1)L) \mbox{ and } \cAt= (1,0)K\vee \infty.
   \end{displaymath}
   Then the mapping (\ref{def-rho}) becomes
   \begin{equation}
   \Pk((l_0,l_1)L)\mapsto (1,l_1 l_0^{-1})K.
   \end{equation}
We shall identify $\cA$ with $L$ via%
   \footnote
   {This is accordance with the inhomogeneous notation used in \cite
   {Benz1}.}
$(1,l)K\equiv\l$. Thus $L$ gets the structure of an affine plane over
$K$. We shall emphasize this by writing $\AkL$ rather than $L$.

   \begin{theo} \label{theo-a}
   Let $\kappa$ be an automorphic collineation or duality of $\Slk$
   fixing $\infty$. Then there exist elements $m_0,m_1,m \in L$,
   $m_0,m_1\neq 0$ and an automorphism or antiautomorphism $J$ of $L$
   with $K^J=K$ such that
      \begin{equation}\label{aff-trafo}
      x^{\rho^{-1} \kappa \rho} = m_1 x^J m_0 + m
      \quad \mbox{for all } x\in L.
      \end{equation}
   The additional conditions
      \begin{eqnarray}
      \label{aff-1}
      {} & J \mbox{ is an automorphism of } L, & {} \\
      \label{aff-2}
      {} & m_0\in K \mbox{ or, only if }K/Z \mbox{ is Galois, }
      m_0i^{-1} \in K & {}
      \end{eqnarray}
   together are necessary and sufficient for $\rho^{-1} \kappa \rho$
   to be an affinity of $\AkL$.
   \end{theo}
{\it Proof.} The assertion in formula (\ref{aff-trafo}) is obviously
true.

Now suppose that $\rho^{-1} \kappa \rho$ is an affinity of $\AkL$.
Then $\kappa$ has to take each chain with a transversal line in $\cAt$
to a chain with a transversal line in $\cAt$. Hence $\cAt^\kappa =
\cAt$ or, only if $K/Z$ is Galois, $\cAt^\kappa = \cAt^{\iota}$.
Therefore $\kappa$ cannot be a duality, so that $J$ cannot be an
antiautomorphism \cite[Theorem 4]{HH3}. Consequently, $g:x \mapsto m_1
x^J m_0$ has to be a semilinear mapping of the right vector space $L$
over $K$. We infer from
   \begin{displaymath}
   xk \stackrel{g}{\mapsto}  (m_1 x^J m_0) (m_0^{-1} k^J m_0)
   \quad \mbox{for all } x \in L, \mbox{ } k \in K
   \end{displaymath}
that $m_0^{-1} K m_0=K$. There are two possibilities: If
   \begin{displaymath}
   m_0^{-1}km_0 = k \quad \mbox{for all } k\in K
   \end{displaymath}
then $m_0$ is a non--zero element of $K$, since $K$ is a maximal
commutative subfield of $L$. On the other hand, however only if $K/Z$
is Galois, also
   \begin{displaymath}
   m_0^{-1}km_0 = \overline{k} \quad \mbox{for all } k\in K
   \end{displaymath}
is possible. Now, again using that $K$ is maximal commutative, it
follows from (\ref{galois-comm}) that $m_0i^{-1} \in K$.

The proof of the converse is a straightforward calculation.\beweisende

\subsection{}
If $\cC$ is a chain such that $(\cC\setminus\{\infty\})^\rho$ is not a
line of $\AkL$ then $(\cC\setminus\{\infty\})^\rho$ will be named a
{\bf circle}. There are two kinds of circles: If $\infty\in \cC$ then
the circle is called {\bf degenerate}, otherwise {\bf
non--degenerate}. The following Lemma shows that distinct chains
cannot define the same circle. In addition it establishes that a
circle cannot be degenerate and non--degenerate at the same time:

   \begin{lemma}\label{lemma-a}
   Let $\cC_0$ and $\cC_1$ be two chains such that
   $\cC_0\setminus\{\infty\} =\cC_1\setminus\{\infty\}$. Then
   $\cC_0=\cC_1$.
   \end{lemma}
{\it Proof}. According to (\ref{sphere-1}), (\ref{sphere-2}),
(\ref{sphere-3}) there exists a 3--dimensional subspace ${\cX}_0$ of
$\Pkh$ with
   \begin{displaymath}
   \cC_0{}^{\ga}={\cX}_0\cap \Pi_Z\cap {\cQ}.
   \end{displaymath}
Since $\cC_0{}^{\ga}$ is an oval quadric of ${\cX}_0\cap \Pi_Z$ and
$Z$ is infinite, $(\cC_0\setminus\{\infty\})^{\ga}$ is still spanning
${\cX}_0$. Repeating this, mutatis mutandis, for $\cC_1$ gives
${\cX}_0={\cX}_1$, whence $\cC_0=\cC_1$, as required.\beweisende

\subsection{}
By Lemma \ref{lemma-a}, we may unambiguously speak of a line being
{\bf tangent} to a circle at some point $P\in\AkL$ or of circles {\bf
touching} at $P$ if they arise from chains that are tangent at
$P^{\rho^{-1}}$.

A degenerate circle has no tangent lines. A point $P$ of a
non--degenerate circle is called {\bf regular} if there exists a
tangent line of that circle at $P$. If such a circle is given as
$\cC^\rho$, $\cC$ a chain, then $P\in \cC^\rho$ is regular if, and
only if, $P$ (regarded as point of ${\cA}$) is incident with a
transversal line of $\cC$. Thus a non--degenerate circle has either
one or two regular points.

\subsection{}
If $K/Z$ is Galois then call two lines, or a circle and a line, or two
circles of $\AkL$ {\bf orthogonal} if they arise from orthogonal
chains.

By virtue of the collineation $\iota$ (cf. formula (\ref{def-iota})),
a line $lK+m$ ($l,m\in L$, $l\neq 0$) is orthogonal to all lines being
parallel to $liK$.

We introduce a unitary scalar product $\ast$ on the right vector space
$L$ over $K$ by setting
   \begin{equation}
(u+iv)\ast (u'+iv')  := \overline{u}u'+\mu_2\overline{v}v'
\quad\mbox{ for all } u,u',v,v'\in K.
   \end{equation}
This scalar product is describing the orthogonality relation on lines
from above. Moreover, $(u+iv)\ast (u+iv)= N(u+iv)$, whence the norm is
a Hermitian form%
   \footnote
   {If $K/Z$ is not Galois then the norm does not seem to be a
   quadratic or Hermitian form on $L$ over $K$.}
on $L$ over $K$.

It is easily seen that there exists no line orthogonal to a degenerate
circle. The join of the two regular points of a non--degenerate circle
is the only line being orthogonal to that circle. It will be called
the {\bf midline} of the circle. The midline is orthogonal to both
tangent lines.

All affinities described in Theorem \ref{theo-a} are preserving
orthogonality.

   \subsection{}
Let $\cC$ be a chain such that $\De:=(\cC\setminus \{ \infty
\})^{\rho}$ is a degenerate circle. Then either there are two points
or there is one point on the line $\infty$ incident with transversal
lines of $\cC$. We call these points at infinity of $\AkL$ the {\bf
absolute points} or the {\bf absolute directions} of $\De$. This
terminology will be motivated in \ref{absol}.

The group $\AGL$ of all transformations (\ref{aff-trafo}) with $m_0=1$
operates sharply 2-fold transitively on $\AkL$. Thus each degenerate
circle can be transferred under $\AGL$ to a degenerate circle through
0 and 1. Write
   \begin{displaymath}
   L^\circ :=
   \left\{
      \begin{array}{l@{\mbox{ : }}l}
      L \setminus (K \cup Ki)       & K/Z\mbox{ Galois,}\\
      L^\circ := L \setminus K      & K/Z\mbox{ not Galois.}
      \end{array}
   \right.
   \end{displaymath}
Then, by \cite[p.329]{Benz1} and (\ref{aff-2}), the degenerate circles
through 0 and 1 are exactly the sets
   \begin{equation}\label{d-c}
   cKc^{-1}\quad\mbox{with }c\in L^\circ .
   \end{equation}

From now on assume that a degenerate circle $\De$ is given by
(\ref{d-c}). Let $\cC$ be the chain with transversal line
$(c,0)K\vee(0,c)K$. Then $\De = (\cC\setminus\{\infty\})^\rho$, whence
$cK$ is an absolute direction of $\De$. Each affinity of $\AGL$ (cf.
formula (\ref{aff-trafo})) with $m_1,m\in cKc^{-1}$ ($m_1\neq 0$,
$m_0=1$ as before) takes $\De$ onto $\De$.

   \begin{theo}
   Each degenerate circle of $\AkL$ is an affine Baer subplane of
   $\AkL$ with the centre of $L$ as underlying field.
   \end{theo}
{\it Proof}. It is sufficient to show this for a degenerate circle
given by (\ref{d-c}). Set $c^{-1}=:d+ie$ with $d,e\in K$. Then, by
(\ref{galois-comm}) and (\ref{non-galois-comm}),
   \begin{displaymath}
   cKc^{-1} =
   \left\{
      \begin{array}{l@{\mbox{ : }}l}
      \{ (cd)k + (cie)\overline{k} \mid k\in K\}   & K/Z \mbox{
      Galois,}\\
      \{ k + (ce)k^D \mid k\in K\}                 &K/Z \mbox{ not
      Galois.}
      \end{array}
   \right.
     \end{displaymath}
Now the assertion follows by Lemma \ref{lemma-b}.\beweisende

\subsection{}
Next we turn to non--degenerate circles.

   \begin{theo}\label{orbit}
All non--degenerate circles of the affine plane $\AkL$ are in one
orbit of $\AGL$.
   \end{theo}
{\it Proof}. Let $\cC_0$ be the chain with transversal line
   \begin{equation}\label{1-i}
   (1,0)K\vee (i,i)K.
   \end{equation}
Then $\Ga_0$:=$\cC_0{}^\rho$ is a non--degenerate circle with regular
point 0.

Let $K/Z$ be Galois. Then 1 is the other regular point of $\Ga_0$. If
$\Ga_1$ is a non--degenerate circle then there exists an affinity
$\al\in\AGL$ taking the regular points of $\Ga_1$ to 0 and 1,
respectively. Hence $\Ga_1{}^{\al\rho^{-1}}$ is a chain with one
transversal line through $(1,0)K$ and the other transversal line
through $(1,1)K$. Applying the collineation $\iota$ on $(1,1)K$
establishes that (\ref{1-i}) is a transversal line of this chain,
whence $\Ga_0=\Ga_1{}^\al$.

Now assume that $K/Z$ is not Galois. If $\Ga_1$ is a non--degenerate
circle then there exists an affinity $\al\in\AGL$ taking the only
regular point of $\Ga_1$ to 0. The chain $\Ga_1{}^{\al\rho^{-1}}$ has
a unique transversal line through $(1,0)K$ and some point of the plane
$(i,0)K\vee \infty$, say
   \begin{displaymath}
   (id,e+if)K \quad\mbox{with } d,e,f\in K,\mbox{ } d,e+if\neq 0.
   \end{displaymath}
There exists an element $m_1\in L\setminus \{ 0\}$ such that
$m_1(e+if) = i d$. The collineation $\kappa$ of $\Pk$ given by
$(l_0,l_1)K \mapsto (l_0,m_1 l_1)K$ leaves $\Slk$ invariant, fixes the
point $(1,0)K$ as well as the line $\infty$ and takes $(id,e+if)K$ to
$(i,i)K$. Hence the induced affinity $\rho^{-1}\kappa\rho$ of $\AkL$
carries $\Ga_1{}^\al$ over to $\Ga_0$.\beweisende

\subsection{}
The non--degenerate circle $\Ga_0$ arising from the chain $\cC_0$ with
transversal line (\ref{1-i}) has the parametric representation
   \begin{equation}\label{n-c}
   \{ ik_1(k_0 + i  k_1)^{-1}\mid (0,0)\neq (k_0,k_1)\in K^2 \};
   \end{equation}
cf. also \cite[Satz 3.2]{Benz1}. Next we establish an equation for
$\Ga_0$:

   \begin{theo}\label{gleich-theo}
   The non--degenerate circle $\Ga_0$ given by (\ref{n-c}) equals the
   set of all points $u+iv$ ($u,v\in K$) satisfying\/%
   \footnote
   {In the elementary plane of complex numbers the same kind of
   equation gives a circle through 0 and 1.}%
      \begin{equation}\label{gleichung}
      u = N(u+iv).
      \end{equation}
   \end{theo}
{\it Proof.} The term $ik_1(k_0+ik_1)^{-1}$ in formula (\ref{n-c}) can
be rewritten as follows: If $K/Z$ is Galois then
   \begin{eqnarray*}
   ik_1(k_0+ik_1)^{-1}      & = &
   ik_1(\overline{k_0}-ik_1)\left((k_0+ik_1)(\overline{k_0}-ik_1)
   \right)^{-1}\\
   {}                       & = &
   (\mu_2k_1\overline{k_1}+i\overline{k_0}k_1)
   (k_0\overline{k_0}+\mu_2k_1\overline{k_1}) ^{-1},
   \end{eqnarray*}
otherwise
   \begin{eqnarray*}
   ik_1(k_0+ik_1)^{-1}      & = &
   ik_1(k_0+k_1+k_1i)\left((k_0+ik_1)(k_0+k_1+k_1i)\right)^{-1}\\
   {}                       & = &
   (\mu_2k_1^2+ik_0k_1)
   \left(k_0^2+k_0k_1+(k_0k_1)^D+\mu_2k_1^2\right) ^{-1}.
   \end{eqnarray*}
Now, since
   \begin{displaymath}
   N(u+iv) =
   \left\{
      \begin{array}{l@{\mbox{ : }}l}
      u\overline{u} + \mu_2 v\overline{v}
            & K/Z\mbox{ Galois,}\\
      u^2 + uv + (uv)^D + \mu_2v^2
            & K/Z\mbox{ not Galois,}
      \end{array}
   \right.
   \end{displaymath}
it is easily seen that all points of $\Ga_0$ are satisfying equation
(\ref{gleichung}).

Conversely, let $q+ir$ ($q,r\in K$) be a solution of
(\ref{gleichung}). If $q=0$ then $r=0$, whence we have a point of
$\Ga_0$. Otherwise set
   \begin{displaymath}
   k_0:=
   \left\{
      \begin{array}{l@{\mbox{ : }}l}
      \mu_2 \overline{rq^{-1}}       & K/Z\mbox{ Galois,} \\
      \mu_2 rq^{-1}                  & K/Z\mbox{ not Galois,}
      \end{array}
   \right\}
   \mbox{ and } k_1:= 1.
   \end{displaymath}
The point of $\Ga_0$ with these parameters equals $q+ir$.\beweisende

\subsection{}
We are able to say a little bit more about non--degenerate circles
provided that $K/Z$ is Galois. Formula (\ref{gleichung}) becomes
   \begin{equation}
   N(u+iv) - u = (u-1+iv)\ast(u+iv) = 0.
   \end{equation}
Thus, if we intersect each line through 0 with its orthogonal line
through 1 then the set of all such points of intersection equals
$\Ga_0$. This is a nice analogon to a well--known property of opposite
points on a Euclidean circle%
   \footnote
   {The points 0 and 1 are, however, the only points of $\Ga_0$ with
   this property.}%
.

   \begin{theo}
   Let $K/Z$ be Galois. Write $E :=\{y\in K \mid y+\overline{y}=1\}$
   and $\cH_e$ ($e\in E$) for the affine Hermitian variety formed by
   all points $u+iv$ ($u,v\in K$) subject to the equation
      \begin{displaymath}
      N(u+iv) = eu+ \overline{eu}.
      \end{displaymath}
   Then the non--degenerate circle $\Ga_0$ given by (\ref{n-c}) can be
   written as
      \begin{equation}
      \Ga_0 = \cH_e\cap\cH_f \quad\mbox{for all }e,f\in E\mbox{ with }
      e\neq f.
      \end{equation}
   \end{theo}
{\it Proof}. A straightforward calculation yields
   \begin{displaymath}
   E = \left\{
      \begin{array}{l@{\mbox{ : }}l}
      \frac{1}{2} + (\la_1+2a)Z     & \mbox{Char}K\neq 2,\\
      a\la_1^{-1} + Z               & \mbox{Char}K   = 2,
      \end{array}
   \right.
    \end{displaymath}
whence $E$ is infinite. Given $q+ir\in\Ga_0$ ($q,r\in K$) then $q\in
Z$ implies
   \begin{displaymath}
   \Ga_0\subset\bigcap_{e\in E} \cH_e.
   \end{displaymath}
Choose distinct elements $e,f\in E$ and $q+ir\in \cH_e\cap\cH_f$
($q,r\in K$). Then
   \begin{displaymath}
   N(q+ir) - N(q+ir) =
   e q+ \overline{e q} - fq- \overline{fq} = 0.
   \end{displaymath}
But
   \begin{displaymath}
   \frac{e-f}{\overline{f}-\overline{e}}=1,
   \end{displaymath}
so that $q=\overline{q}$ and therefore $q+ir\in \Ga_0$.\beweisende

\subsection{}
There is an alternative approach to $\AkL$ via the point model of
$\Sigma(K,L)$ on the Klein quadric $\cQ$.

Write $I:=\infty^\ga$ and $\cZ$ for the $\ga$--image of the ruled
plane on $\cAt$; this $\cZ$ is a plane on the Klein quadric.
Furthermore let $\cFt$ be any plane of $\Pkh$ skew to $\cZ$ and write
   \begin{equation}
   \pi : \Pkh\setminus\cZ \to \cFt
   \end{equation}
for the projection with centre $\cZ$ onto the plane $\cFt$. It is well
known from descriptive line geometry that there exists a collineation
$\psi $ of $\cAt$ onto $\cFt$ such that
   \begin{displaymath}
   (p\cap\cAt)^\psi = p^{\ga\pi}
   \end{displaymath}
for all lines $p$ of $\Pk$ not contained in $\cAt$. Cf., e.g.,
\cite{Brau1}. We turn $\cFt$ into an affine plane $\cF$, say, by
regarding $\cFt\cap I^\perp$ as its line at infinity; here $I^\perp$
denotes the tangent hyperplane of the Klein quadric at $I$. Then
$\infty^\psi =\cF\cap I^\perp$.

The bijectivity of $\rho$ implies that $\Slk{}^\ga\setminus\{ I \}$ is
mapped bijectively under $\pi$ onto the affine plane $\cF$. The
restriction
   \begin{displaymath}
   \pi\mid\Slk{}^\ga\setminus\{ I \}
   \end{displaymath}
can be seen as a {\bf generalized stereographic projection} of the
oval quadric $\Slk{}^\ga$ of $\Pi_Z$ onto the affine plane%
   \footnote
   {A `usual' stereographic projection would map onto a 4--dimensional
   affine space over $Z$ rather than an affine plane over $K$.}
$\cF$.

Let $\cC$ be a chain. Then $\cC^\ga = \cX\cap{\cQ} \cap \Pi_Z$ for
some 3-dimensional subspace $\cX$ of $\Pkh$. We leave it to the reader
to show that $(\cC\setminus\{\infty\})^{\ga\pi}$ is an affine line if
$\cX\cap\cZ$ is a line through $I$, a degenerate circle if $\cX\cap\cZ
= \{I\}$ and a non--degenerate circle if $\cX\cap\cZ$ is some point
other than $I$.

Using the mapping $\ga\pi\psi^{-1}$ instead of $\rho$ is very
convenient to establish results on the images of traces
\cite[p.327]{Benz1}, since their $\ga$--images are just the regular
conics on $\Slk{}^\ga$ \cite[3.4]{HH3}. We sketch just one result
without proof:

Let $\cC$ be a chain through $\infty$ such that $(\cC\setminus\{
\infty \})^\rho =: \De$ is a degenerate circle of $\AkL$. Then the
$\rho$--images of traces in $\cC$ are on one hand the lines of the
affine plane $\De$ and on the other hand certain ellipses of $\De$. If
these ellipses are extended to conics of $\AkL$ then the absolute
directions of $\De$ determine their points at infinity\label{absol}%
   \footnote
   {There is only one such point if $K/Z$ is not Galois.}%
. This is the well--known concept of absolute circular points. $\De$
is a Euclidean plane representing the extension $K/Z$. Cf.
\cite{Schroe1}.

%%%%%%%%%%%%%%%%%%%%%%%%%%%%%%%%%%%%%%%%%%%%%%%%%%%%%%%%%%%%%%%%%%%%%%
% Hyphenation of some German words is given explicitely below, since
% English patterns yield wrong results.
%%%%%%%%%%%%%%%%%%%%%%%%%%%%%%%%%%%%%%%%%%%%%%%%%%%%%%%%%%%%%%%%%%%%%%
   
   \small
Hans Havlicek, Abteilung f\"ur Lineare Algebra und Geometrie,
Technische Universit\"at, Wiedner Hauptstrasse 8--10, A--1040 Wien,
Austria.

   \noindent %%%%%%%%%%%%%%%%%%%%%%%%%%%%%%%%%%%%%%%%%%%%%%%%%%%%%%
EMAIL: havlicek@geometrie.tuwien.ac.at

\end{document}